\documentclass[reqno,12pt]{amsart}
\usepackage{amsmath,amsthm,amssymb,tabularx}

\theoremstyle{plain}

\newtheorem{tm}{Theorem}[section]
\newtheorem{lm}[tm]{Lemma}
\newtheorem{cor}[tm]{Corollary}
\theoremstyle{definition}
\newtheorem{ex}[tm]{Example}
\newtheorem{prop}[tm]{Proposition}
\newtheorem{df}[tm]{Definition}

\newcommand{\beq}{\begin{equation}}
\newcommand{\eeq}{\end{equation}}
\newcommand{\bga}{\begin{gather*}}
\newcommand{\ega}{\end{gather*}}
\newcommand{\bit}{\begin{itemize}}
\newcommand{\eit}{\end{itemize}}
\newcommand{\btm}{\begin{tm}}
\newcommand{\etm}{\end{tm}}
\newcommand{\blm}{\begin{lm}}
\newcommand{\elm}{\end{lm}}
\newcommand{\bcor}{\begin{cor}}
\newcommand{\ecor}{\end{cor}}
\newcommand{\bex}{\begin{ex}}
\newcommand{\eex}{\end{ex}}
\newcommand{\bpr}{\begin{proof}}
\newcommand{\epr}{\end{proof}}
\newcommand{\bprop}{\begin{prop}}
\newcommand{\eprop}{\end{prop}}
\newcommand{\bdf}{\begin{df}}
\newcommand{\edf}{\end{df}}

\def\R{\mathbb{R}}

\def\L{{\mathcal L}}

\def\Z{\mathbb{Z}}

\let\a\alpha
\let\b\beta

\let\g\gamma
\def\e{\varepsilon}
\let\s\sigma

\let\wt\widetilde

\def \le {\leqslant}
\def \ge {\geqslant}
\let \<\langle
\let \>\rangle

\let\phi\varphi
\let\trg\triangledown
\def\lpw{\L_p^w(G)}
\def\lqw{\L_q^{w^{\!-\!1}}(G)}
\def\qw{{q,w^{-\!1}}}
\def\pw{{p,w}}
\def\sk#1{\noalign{\kern#1}}

\begin{document}

\title[Example of a weighted algebra on uncountable group]
{AN EXAMPLE OF A WEIGHTED ALGEBRA $\L_2^w(G)$ ON UNCOUNTABLE DISCRETE GROUP}
 \author{Yu. N. Kuznetsova}
 \address{VINITI, Mathematics dept., Usievicha 20, Moscow 125190}
 \email{jkuzn@mccme.ru}
\thanks{Work supported by RFBR grant No. 08-01-00867.}

\maketitle

Let $G$ be a locally compact group, and let $w>0$ be a measurable function.
We consider weighted spaces $\lpw$ with $p\ge1$, which are by definition equal
$$
\lpw=\{f: \int_G |fw|^p<\infty\}.
$$
Norm of a function $f$ in this space is $\|f\|_\pw=\|fw\|_p=(\int_G|fw|^p)^{1/p}$.

Sufficient conditions are well-known, under which a weighted space becomes an
algebra with respect to convolution. For $p=1$ it is
sub\-mul\-ti\-pli\-ca\-ti\-vi\-ty \cite{grsh}:
\beq\label{eq_submult}
w(st) \le w(s)w(t),
\eeq
and for $p>1$ the following inequality \cite{wermer} (pointwise locally almost everywhere):
\beq\label{eqw}
w^{-q} * w^{-q} \le w^{-q}.
\eeq
Here $q$ is the conjugate exponent to $p$, so that $1/p+1/q=1$.
Commu\-ta\-ti\-vi\-ty of such an algebra is equivalent to that of the group.

Natural example of a weighted algebra is given by the usual group algebra
$\L_1(G)$ with trivial weight $w\equiv1$. But if $p>1$, existence
of weighted algebras on this or that group is a nontrivial question.
For example, the usual spaces $\L_p(G)$ are closed under convolution on compact
groups only (this is the positive solution of renowned $L_p$-hypothesis,
see final result in Saeki \cite{saeki}).

We are interested in possibly full description of the class $\mathcal {WL}_p$
of groups on which weighted algebras $\lpw$ with $p>1$ exist. It has long been
known that the real line belongs to this class \cite{wermer} (one may take,
e.g., $w(t)=1+t^2$ for any $p>1$), as well as $\Z$ and $\R^n$.
In \cite{mz} it is proved that this class contains all $\sigma$-compact groups
(i.e. representable as a countable union of compact sets).
There is also proved that for abelian groups $\s$-compactness is necessary
for containment in this class.

In the present paper we show that $\mathcal {WL}_p$ for $1<p\le 2$ is wider than
the class of $\s$-compact groups. Theorem \ref{free} presents a construction
of weighted algebras on an uncountable discrete group (in the discrete case
$\s$-compactness is the same as countability). For $p>2$ there are no weighted
algebras on these groups, so that the classes $\mathcal {WL}_p$ are different
for different $p$.

\S2 treats weighted algebras with involution. The results obtained allow to
extend the theorem claiming that $\sigma$-compactness is necessary for the
existence of weighted algebras, to the class of amenable groups. This is done in \S3.
Of course, groups in theorem \ref{free} are not amenable; they are free groups
with uncountable set of generators.

The results obtained give answer also to the following question. It was known
before \cite{kuz-faa} that inequality \eqref{eqw} may not hold for a weight
of an algebra with $p>1$. But in the examples of \cite{kuz-faa} algebras are
not translation invariant, and weights are not locally bounded. Do these
additional assumptions imply \eqref{eqw}? Theorem \ref{free} answers this
question in negative.

\section{Example of a weighted algebra on a free group}

\btm\label{free}
On a free group $F$ of any infinite cardinality there is a weight $w$
such that weighted space $\L_p^w(F)$, $1<p\le 2$, is a Banach algebra with respect
to convolution.
\etm

\bpr
Denote by $A$ the set of generators of $F$. We consider elements of $F$
as reduced words in the alphabet $A\cup A^{-1}$. Let $|\a|$ denote
the length of a word $\a\in F$, and put $A_n=\{\a:|\a|=n\}$, $n\ge0$. We define
$$
w|_{A_n}=(n+1)^3
$$
(any number $>2$ may be taken instead of 3).

In this example we denote $\ell_p=\L_p(F)$, $\ell_p^w=\L_p^w(F)$.
Since $\ell_p^w=\ell_p/w$, it is sufficient to take $f,g\in\ell_p$
and show that
$$
h=w\cdot\left({f\over w}*{g\over w}\right)\in\ell_p.
$$
We may assume that $f,g$ are nonnegative. For any $\a$
$$
h(\a)=w(\a)\cdot\left({f\over w}*{g\over w}\right)(\a)=
 w(\a)\sum_\b{f_\b g_{\b^{-1}\a}\over w(\b)w(\b^{-1}\a)}.
$$
Let $\a\in A_n$, $\b\in A_k$.
\bit
\item[1)] If $k\ge n/2$, then $w(\b)\ge (n/2+1)^3\ge (n+1)^3/8=w(\a)/8$,
  and $w(\a)/w(\b)\le8$.
\item[2)] If $k< n/2$, then $\b^{-1}\a$ contains no less than $n-k$ letters,
so that $w(\b^{-1}\a)\ge(n-k+1)^3>(n/2+1)^3\ge w(\a)/8$, and
then $w(\a)/w(\b^{-1}\a)\le8$.
\eit
Thus, we may write an estimate $h(\a)\le8\phi(\a)+8\psi(\a)$, where
\begin{gather*}
\phi(\a)=\sum_\b {f_\b g_{\b^{-1}\a}\over w(\b)}, \\
\psi(\a)=\sum_\b {f_\b g_{\b^{-1}\a}\over w(\b^{-1}\a)}
 = \sum_\b {f_{\a\g^{-1}} g_\g\over w(\g)}.
\end{gather*}
One may see later that $\phi$ and $\psi$ are estimated
similarly, and it will suffice to estimate $\phi$ alone.
\begin{gather*}
\phi(\a)=\sum_\b {f_\b g_{\b^{-1}\a}\over w(\b)}=
 \sum_{k=0}^\infty \sum_{\b\in A_k} {f_\b g_{\b^{-1}\a}\over (k+1)^3}.
\end{gather*}
Each set $A_k$ may be split into disjoint sets
$A_{kj}(\a)$, $0\le j\le k$, by number of common letters in $\a$ and $\b$:
$$
A_{kj}(\a)=\{\b\in A_k: \b_1=\a_1,\ldots,\b_j=\a_j,\b_{j+1}\ne\a_{j+1}\}.
$$
If $j>n=|\a|$, we assume that $A_{kj}(\a)=\emptyset$.
Thus, $A_k=\cup_{j=0}^k A_{kj}(\a)$, so  that
$$
\phi(\a)= \sum_{k=0}^\infty \sum_{j=0}^k
 \sum_{\b\in A_{kj}(\a)} {f_\b g_{\b^{-1}\a}\over (k+1)^3}
 = \sum_{k=0}^\infty\sum_{j=0}^k {\phi_{kj}(\a)\over (k+1)^3},
$$
where
$$
\phi_{kj}(\a)=\sum_{\b\in A_{kj}(\a)} f_\b g_{\b^{-1}\a}.
$$
In other words,
$$
\phi=\sum_{k=0}^\infty\sum_{j=0}^k {\phi_{kj}\over (k+1)^3}.
$$
Estimate now the norms $\|\phi_{kj}\|_p$, taking into account that $\phi_{kj}(\a)=0$
when $n=|\a|<j$. Now, and only now, we use the fact that $p\le2$: it implies
that $\|x\|_q\le\|x\|_p$ for the conjugate exponent $q$ and all $x\in\ell_p$.
\begin{gather*}
\|\phi_{kj}\|_p^p =\sum_{n=0}^\infty \sum_{\a\in A_n} \phi_{kj}(\a)^p
= \sum_{n=j}^\infty \sum_{\a\in A_n} \phi_{kj}(\a)^p
=\\=
\sum_{n=j}^\infty \sum_{\a\in A_n} \bigg(\sum_{\b\in A_{kj}(\a)} f_\b g_{\b^{-1}\a}\bigg)^p
 \le \sum_{n=j}^\infty \sum_{\a\in A_n}
  \sum_{\b\in A_{kj}(\a)} f_\b^p \sum_{\g\in A_{kj}(\a)} g_{\g^{-1}\a}^p.
\end{gather*}
A word $\a\in A_n$ may be represented in the form $\a=\hat\a\check\a$, where
$\hat\a=\a_1\ldots\a_j\in A_j$, $\check\a=\a_{j+1}\ldots\a_n\in A_{n-\!j}$
and $\a_{j+1}\ne\a_j^{-1}$.
Then every $\b\in A_{kj}(\a)$ has form $\b=\hat\a\hat\b$, where
$\hat\b=\b_{j+1}\ldots\b_k\in A_{k-j}$ and $\b_{j+1}\ne\a_j^{-1}$, $\b_{j+1}\ne\a_{j+1}$.
Similar representation holds for $\g=\hat\a\hat\g$,
so that $\g^{-1}\a=\hat\g^{-1}\check\a$.
Thus,
\begin{gather*}
\|\phi_{kj}\|_p^p
\le \sum_{n=j}^\infty \sum_{\hat\a\in A_j}
 \sum_{\scriptstyle \check\a\in A_{n-\!j}\atop
  \scriptstyle \a_{j+1}\ne\a_j^{-1}}
 \sum_{\scriptstyle\hat\b\in A_{k-\!j}\atop
  \scriptstyle \b_{j+1}\ne\a_j^{-1}\!,\, \b_{j+1}\ne\a_{j+1}}
 \mkern-18mu  f_{\hat\a \hat\b}^p
 \sum_{\scriptstyle \hat\g\in A_{k-\!j}\atop
  \scriptstyle \g_{j+1}\ne\a_{j+1}\!,\, \g_{j+1}\ne\a_j^{-1}}
 \mkern-18mu g_{\hat\g^{-1}\check\a}^p
  \le\\
  \le \sum_{n=j}^\infty \sum_{\hat\a\in A_j}
 \sum_{\scriptstyle\hat\b\in A_{k-\!j}\atop \scriptstyle \b_{j+1}\ne\a_j^{-1}}
  f_{\hat\a\hat\b}^p \sum_{\check\a\in A_{n-\!j}}
 \sum_{\scriptstyle \hat\g\in A_{k-\!j}\atop \scriptstyle \g_{j+1}\ne\a_{j+1}}
  g_{\hat\g^{-1}\check\a}^p
  =\\
  =\sum_{n=j}^\infty\,\sum_{\xi=\hat\a\hat\b\in A_k} f_\xi^p
   \sum_{\eta=\hat\g^{\!-\!1}\!\check\a\in A_{n+k-2j}} g_\eta^p
    \le\|f\|_p^p \,\sum_{n=j}^\infty \,\sum_{\eta\in A_{n+k-2j}} g_\eta^p
  \le\|f\|_p^p\cdot\|g\|_p^p.
\end{gather*}
That is, $\|\phi_{kj}\|_p \le \|f\|_p\cdot\|g\|_p$. Hence
\begin{gather*}
\|\phi\|_p\le \sum_{k=0}^\infty\sum_{j=0}^k{\|\phi_{kj}\|_p\over (k+1)^3}
\le \sum_{k=0}^\infty(k+1){\|f\|_p\cdot\|g\|_p\over (k+1)^3} = \\
=\sum_{k=0}^\infty {1\over (k+1)^2} \|f\|_p\cdot\|g\|_p \equiv C\|f\|_p\cdot\|g\|_p.
\end{gather*}
The norm $\|\psi\|_p$ is estimated similarly, replacing sum over $\b$ with that over
$\g=\b^{-1}\a$. Thus, $\|h\|_p\le 16C\|f\|_p\cdot\|g\|_p$,
what was to show.
\epr

It is easy to show that examples of this sort for $p>2$ are impossible:

\bprop
If $\lpw$ is an algebra and $G$ is a discrete uncountable group, then $p\le2$.
\eprop

\bpr
Since $G$ is uncountable, for some $C>0$ the set
$A=\{x: \max(w(x),w(x^{-1})\le C\}$ is also uncountable. Note that $A=A^{-1}$.
Now, if $f\in\ell_p(A)$, then
$$
\|f\|_\pw^p =\sum_{x\in A} |f(x)w(x)|^p \le C\sum_{x\in A}\|f\|_p^p,
$$
so that $\ell_p(A)\subset \lpw$. In particular,
$$
\sum_{x\in A} f(x)g(x^{-1}) = (f*g)(e)<\infty
$$
for all nonnegative $f,g\in\ell_p(A)$. We have now $\ell_p(A)\cdot\ell_p(A)\subset\ell_1(A)$,
what is possible for $p\le2$ only.
\epr

\section{Symmetric weighted algebras}\label{symmetr}

Natural involution is not defined in all algebras $\lpw$,
but in certain ones which we call symmetric. After series of lemmas
we prove that symmetric algebras posess an injective involutive
representation in $\L_2(G)$, and are as a consequence semisimple.
That commutative, not necessarirly symmetric algebras are semisimple,
was proved in \cite[th. 4]{kuz-faa}.

In this section $\Delta$ denotes modular function of $G$.
We use also notation ${}^\trg$: $\,f^\trg(t)=f(t^{-1})$.

\bdf
Define
\beq\label{invol}
f^*(t)=\bar f(t^{-1})\Delta(t^{-1}).
\eeq
A space $\L_p^w(G)$ is called {\it symmetric}, if $f^*\in \L_p^w(G)$
and $\|f\|_\pw=\|f^*\|_\pw$ for all $f\in \L_p^w(G)$.
\edf

The set of all functions of the type $f^*$, where $f\in\lpw$,
constitutes a weighted space $\L_p^{\wt w}(G)$
with the weight $\wt w=\Delta^{(p-1)/p}w^\trg$,
whereas $\|f^*\|_\pw = \|f\|_{p,\wt w}$:
$$
\|f^*\|^p_\pw = \int |\bar f(t^{\!-1})\Delta(t^{\!-1})w(t)|^pdt
 = \int |f(t)w(t^{\!-1})|^p\Delta(t)^{p-1}dt
 = \|f\|^p_{p,\wt w}.
$$
Sufficient condition for symmetry is equality $w=\wt w$,
which turns to usual evenness in the case of $p=1$ or a unimodular group.

\blm\label{inv_alg}
If the space $\L_p^w(G)$ is an algebra, then the space $\L_p^{\wt w}(G)$
with the weight $\wt w=\Delta^{(p-1)/p}w^\trg$ is also an algebra.
\elm

\bpr
This follows from the fact that $f^* *g^*=(g*f)^*$ for all $f,g\in \L_p^w(G)$.
Moreover, obviusoly, $\|\phi*\psi\|_{p,\wt w}\le \|\phi\|_{p,\wt w}\|\psi\|_{p,\wt w}$
for all $\phi,\psi\in \L_p^{\wt w}(G)$.
\epr

\blm\label{inv_max}
If the spaces $\L_p^w(G)$, $\L_p^v(G)$ are both algebras, then the space
$\L_p^u(G)$ with the weight $u=\max\{w,v\}$ is also an algebra.
\elm

\bpr
This follows from the fact that $\L_p^u(G) = \L_p^w(G)\cap \L_p^v(G)$.
\epr

\blm\label{sym_exist}
If the space $\L_p^w(G)$ is an algebra, there is a weight $v$ on $G$
such that the space $\L_p^v(G)$ is a symmetric algebra.
\elm

\bpr
According to lemmas \ref{inv_alg}, \ref{inv_max}, the proper weight is
$$v=\max\{w, \Delta^{-(p-1)/p}w^\trg\}.$$
\epr

\blm\label{sym_module}
If $\L_p^w(G)$ is a symmetric algebra, then $\|f*g\|_\qw\le \|f\|_\pw\cdot\|g\|_\qw$
for all $f\in\lpw$, $g\in\lqw$.
\elm

\bpr
Take arbitrary $f,\phi\in\lpw$, $g\in\lqw$ and transform the following expression:
\begin{align*}
\<\phi, f*g\>
&= \int_G\int_G \phi(t)f(s)g(s^{-1}t)ds\,dt 
= \int_G\int_G \phi(su)f(s)g(u)ds\,du =\\
&= \int_G\int_G \phi(s^{-\!1}u)f(s^{-\!1})\Delta(s^{-\!1})g(u)ds\,du=\\
&= \int_G\int_G \phi(s^{-\!1}u)\bar f^*(s)g(u)ds\,du
= \<\bar f^* *\phi, g\>.
\end{align*}
Since the spaces $\lqw$ are $\lpw$ conjugate to each other, and their natural norms
coincide with those in the sense of conjugate spaces,
$$
\|f*g\|_\qw\! = \sup_{\phi\ne0} {|\<\phi, f*g\>|\over \|\phi\|_\pw}
 = \sup_{\phi\ne0} {|\<\bar f^*\! *\phi, g\>|\over \|\phi\|_\pw}
 \le \sup_{\phi\ne0} {\|\bar f^*\! *\phi\|_\pw \|g\|_\qw \over \|\phi\|_\pw}.
$$
Using the inequality $\|\bar f^*\! *\phi\|_\pw \le \|\bar f^*\|_\pw\cdot\|\phi\|_\pw$,
and symmetricity of $\lpw$, we get:
$$
\|f*g\|_\qw \le \|\bar f^*\|_\pw\|g\|_\qw = \|f\|_\pw \|g\|_\qw,
$$
what was to show.
\epr

Operation $f\mapsto f^*$ is an involution on a symmetric algebra. Similarly to
the case of usual group algebra, this yields a representation in $\L_2(G)$.
In the proof of the next theorem ideas of Kerman and Sawyer \cite{KS} are used.

\btm\label{repres}
If $\L_p^w(G)$ is a symmetric algebra, then:
\bit\item[(1)] Convolution $f*g$ is defined and belongs to $\L_2(G)$
for all $f\in \lpw$, $g\in \L_2(G)$, and $\|f*g\|_2 \le \|f\|_\pw \|g\|_2$;
\item[(2)] Convolution operator $T_f(g)=f*g$ on $\L_2(G)$ is bounded, and $\|T_f\|\le\|f\|_{p,w}$;
\item[(3)] The map $f\mapsto T_f$ is an injective involutive representation of $\lpw$.
\eit
\etm

\bpr
Statement (1) is enough to prove for nonnegative $f$ and $g$.
As any summable function, $g^2$ may be represented in the form $g^2=g_1g_2$,
where $g_1=g^{2/p}\in \L_p(G)$, $g_2=g^{2/q}\in \L_q(G)$. Put $\phi=g_1/w\in\lpw$,
$\psi=g_2w\in\lqw$. Then $g^2=\phi\psi$ and $\|\phi\|_{p,w} = \|\psi\|_{q,w^{-1}}=\|g\|_2$.
Now in every point $s\in G$, using Cauchy-Bunyakovskii inequality, we have:
\begin{gather*}
(f*g)(s) = \int_G f(t) \sqrt{\phi(t^{-1}s)\psi(t^{-1}s)}dt
\le\\
 \le \Big(\int_G f(t)\phi(t^{-1}s)dt \int_G f(r)\psi(r^{-1}s)dr\Big)^{1/2}
 = \sqrt{(f*\phi)(s) (f*\psi)(s)},
\end{gather*}
and, with account of lemma \ref{sym_module},
\begin{gather*}
\|f*g\|_2^2 =
 \int_G |(f*g)(s)|^2ds \le \int_G (f*\phi)(s)(f*\psi)(s) ds \le
\\\le
 \|f*\phi\|_{p,w} \|f*\psi\|_{q,w^{-1}}
 \le \|f\|_{p,w}^2\|\phi\|_{p,w} \|\psi\|_{q,w^{-1}} = \|f\|_{p,w}^2\|g\|_2^2.
\end{gather*}

Thus (1) is proved, and (2) follow trivially.

Proof of (3) is identical to the case of $\L_1(G)$, see, e.g., \cite[\S 28]{naimark}.
\epr

\bcor
Every symmetric algebra $\L_p^w(G)$ is semisimple.
\ecor

\section{Corollaries for amenable groups}

We prove first a weighted analog of the theorem of \.Zelazko \cite{zel}.

\blm\label{muab}
Let $\lpw$ be a symmetric algebra. Then for any sets $A,B$ of positive
finite measure
$$
\mu(AB)\ge \mu(B) \|I_A\|_\qw^2.
$$
\elm

\bpr
Pick any $A$ and $B$ (of finite measure).
Using simple properties of convolution, Cauchy-Bunyakovskii inequality
and theorem \ref{repres}, we get:
\begin{align*}
\|I_A\|_\qw^q \|I_B\|_2^2
&= \int_A w^{-q} \cdot\int_B 1
 = \int_G I_A w^{-q}\cdot\int_G I_B =\\
 &=\int_G(I_Aw^{-q})*I_B=\int_{AB}(I_Aw^{-q})*I_B \le\\
 &\le\|I_{AB}\|_2 \|(I_Aw^{-q})*I_B\|_2 \le
 \|I_{AB}\|_2 \|I_Aw^{-q}\|_\pw\|I_B\|_2=\\
&=\|I_{AB}\|_2 \Big(\int_A w^{-pq+p}\Big)^{1/p} \|I_B\|_2 =\\
&=\|I_{AB}\|_2 \Big(\int_A w^{-q}\Big)^{1/p} \|I_B\|_2
=\|I_{AB}\|_2\|I_A\|_\qw^{q/p} \|I_B\|_2
\end{align*}
(as $-pq+p=-q$). Hence, as $q-q/p=1$,
$$
\|I_A\|_\qw\|I_B\|_2 \le \|I_{AB}\|_2.
$$
Passing to squares, we get the statement of the lemma.
\epr

In the abelian case algebra with an even weight, $w(x)=w(x^{-1})$, is
contained in $\L_1(G)$ (this was proved for $G=\R^n$ by Kerman and Sawyer \cite{KS}).
The next theorem generalizes this result to amenable groups, but the condition
of evenness must be replaced by symmetricity of algebra.
Theorem \ref{free} shows that for arbitrary groups the same result does not hold.

\btm\label{amen_l1}
Let $G$ be an amenable group and let $\lpw$, $p>1$, be a symmetric algebra.
Then $w^{-q}\in\L_1(G)$ and $\lpw\subset \L_1(G)$.
\etm

\bpr
We use the uniform F\o lner condition \cite{em-gr}: for any compact set
$A\subset G$ containing identity, and any $\e>0$ there is a compact set
$B\subset G$ such that $\mu(AB\Delta B)<\e\mu(B)$. Hence it follows that
$\mu(AB)< (1+\e)\mu(B)$. But by lemma \ref{muab} $\mu(AB)\ge \mu(B) \|I_A\|_\qw^2$.
Thus, $\|I_A\|_\qw \le 1$, i.e. ${\int_{\!A} w^{\!-q}\le1}$.
As $A$ was chosen arbitrary, we conclude that $\int_G w^{-q} \le 1$,
i.e. $w^{-q}\in\L_1(G)$. Hence it follows easily \cite[prop. 2]{kuz-faa}
that $\lpw\subset \L_1(G)$.
\epr

Finally, we show that among amenable groups, weighted algebras with $p>1$
exist on $\sigma$-compact groups only.

\bcor
If $G$ is an amenable group and with some $w$ and $p>1$ the space $\lpw$
is an algebra, then the group $G$ is $\sigma$-compact.
\ecor

\bpr
By lemma \ref{sym_exist} we may assume $\lpw$ is symmetric. Then by theorem
\ref{amen_l1} $w^{-q}\in\L_1(G)$. As $w^{-q}$ is strictly positive, and the support
of any summable function may be chosen to be $\sigma$-compact \cite[11.40]{HR},
we get that entire group $G$ is $\sigma$-compact.
\epr

\end{document}